\documentclass[12pt]{article}
\usepackage{amsmath,amsthm,amsfonts,amssymb,amscd}
\usepackage{graphics}
\usepackage[latin1]{inputenc}

\input{tcilatex}
\begin{document}

\renewcommand{\thefootnote}{\fnsymbol{footnote}}

\centerline{\LARGE\bf The cohomology of the weak stable}

\vskip .1in

\centerline{\LARGE\bf foliation of geodesic flows}

\vskip .5in

\centerline{\sc Nathan M. dos Santos}

\vskip .1in

\centerline{\sl Instituto de Matemática, Universidade Federal Fluminense,}

\centerline{\sl 24020-005, Niterói, RJ, Brazil}

\centerline{e-mails: santos@mat.uff.br}

\ 

\section{Introduction}

Let $M=SL(2,R)/{\Gamma }$ where ${\Gamma }$ is a cocompact lattice and $GA$
be the subgroup of upper-triangular matrices. $GA$ is isomorphic of the
group of orientation-preserving affine transformations of the real line. The
Lie algebra $s\ell (2,R)$ of $SL(2,R)$ has the canonical basis 
\begin{equation}
Y=\frac{1}{2}%
\begin{pmatrix}
1 & 0 \\ 
0 & -1%
\end{pmatrix}%
,\quad S=%
\begin{pmatrix}
0 & 1 \\ 
0 & 0%
\end{pmatrix}%
\quad \text{ and }\quad U=%
\begin{pmatrix}
0 & 0 \\ 
1 & 0%
\end{pmatrix}
\tag{1.1}
\end{equation}%
satisfying the structural equations 
\begin{equation}
\lbrack Y,S]=S,\quad \lbrack Y,U]=-U\quad \text{ and }\quad \lbrack S,U]=2Y 
\tag{1.2}
\end{equation}%
and for the dual basis $\{\sigma ,\eta ,\mu \}$ we have 
\begin{equation}
d\sigma =\sigma \wedge \eta ,\quad d\eta =2\mu \wedge \sigma \quad \text{
and }\quad d\mu =\eta \wedge \mu .  \tag{1.3}
\end{equation}%
The right action of $SL(2,R)$ on itself restricts to $GA$ giving the
homogeneous action 
\begin{equation}
A\colon M\times GA\rightarrow M,(s,g)=sg  \tag{1.4}
\end{equation}%
and the orbit foliation ${\mathcal{F}}$ is the \textit{weak stable
foliation\/} of the Anosov flow generated by $Y$. We have the injective
homomorphism 
\begin{equation}
i_{A}\colon \mathcal{\tciLaplace }\rightarrow \chi ^{\infty }({\mathcal{F}})
\tag{1.5}
\end{equation}%
of the Lie algebra $\mathcal{\tciLaplace }$ of $GA$ into the Lie algebra $%
\chi ({\mathcal{F}})$ of the smooth vector fields tangent to the weak stable
foliation ${\mathcal{F}}$ of the geodesic flow ${\varphi }_{t}$ generated by 
$Y$ defined by 
\begin{equation*}
i_{A}(E)(x)=DA_{x}(0)\cdot E=E_{A}(x).
\end{equation*}%
The map $I_{A}$ is defined by 
\begin{equation}
I_{A}\colon {\Lambda }_{R}(\mathcal{\tciLaplace })\rightarrow {\Lambda }({%
\mathcal{F}}),\quad I_{A}(w)\cdot E_{A}=w_{A}(E_{A})=w(E)  \tag{1.6}
\end{equation}%
where ${\Lambda }_{R}(\mathcal{\tciLaplace })$ is the exterior algebra of
the invariant forms of $GA$ into the leafwise forms on $M$. Associated to
the $A$-invariant volume form $\sigma \wedge \eta \wedge \mu $ we have the
chain map $[S]$ 
\begin{equation}
P_{A}:\Lambda (M)\rightarrow {\Lambda }_{R}(\mathcal{\tciLaplace }) 
\tag{1.7}
\end{equation}%
given by 
\begin{equation*}
P_{A}(w)\cdot E=\int_{M}w(E_{A})d\nu ,\quad E\in \mathcal{\tciLaplace }
\end{equation*}%
where $d\nu $ is the measure given by $\sigma \wedge \eta \wedge \mu $ ,and
similarly if $w$ is a 2-form.Since $I(\tciFourier )\subset \ker P_{A}$ \ we
have the chain map \ $P_{A}^{^{0}}:\Lambda (\tciFourier )\rightarrow \Lambda
(\tciLaplace )$

The kernel $K$ of $P_{A}^{o}$ is a chain complex and we have the split short
exact sequence 
\begin{equation}
0\rightarrow K^{j}\,\,\,\overset{i}{\rightarrow }\,\,\,{\Lambda }^{j}({%
\mathcal{F}})\quad 
\begin{matrix}
P_{A}^{o} \\ 
\rightleftarrows \\ 
i_{A}%
\end{matrix}%
\quad {\Lambda }^{j}(\mathcal{\tciLaplace })\rightarrow 0  \tag{1.8}
\end{equation}%
giving 
\begin{equation}
H^{j}({\mathcal{F}})=H^{j}(\mathcal{\tciLaplace })\oplus H^{j}(K),\quad
0<j\leq 2.  \tag{1.9}
\end{equation}%
The cohomology of the Lie algebra $\mathcal{\tciLaplace }$ is given by 
\begin{equation}
H^{1}(\mathcal{\tciLaplace })=R[\eta ]\quad \text{ and }\quad H^{2}(\mathcal{%
\tciLaplace })=0.  \tag{1.10}
\end{equation}

\section{The computation of the cohomology $H^*({\mathcal{F}})$}

The main result of this paper is the computation of $H^2({\mathcal{F}})$.
The computation of $H^1({\mathcal{F}})=R[\eta]\oplus H^1(M,R)$ was done by
S. Matsumoto and Y. Mitsumatsu [M.M].

The \textit{reduced leafwise cohomology\/} ${\mathcal{H}}^{j}({\mathcal{F}})$
of a foliated manifold $(M,{\mathcal{F}})$ is the quotient of the leafwise
closed $j$-forms by the closure of the coboundary$B^{j}({\mathcal{F}})$. J.
Alvares and G. Hector [A.H] have given sufficient conditions for a foliation
to have infinite dimension reduced cohomology and gave plenty of examples.
S. Matsumoto [M] proved that for the weak stable foliation ${\mathcal{F}}$
of a geodesic flow we have 
\begin{equation}
{\mathcal{H}}^{2}({\mathcal{F}})\simeq H^{2}(M,R).  \tag{2.1}
\end{equation}

\vskip .2in

\noindent \textbf{Theorem.} \textit{If ${\mathcal{F}}$ is the weak stable
foliation of a geodesic flow, then 
\begin{equation*}
H^{2}({\mathcal{F}})=H^{2}(M,R).
\end{equation*}%
}

\vskip .2in

\noindent \textit{Proof.} In view of [M., Theorem 13] we have to show that
the coboundary space $B^{2}({\mathcal{F}})$ is closed. We recall the
definition of leafwise forms. Let ${\Lambda }(M)$ be the space of smooth
differential forms on $M$ and $I({\mathcal{F}})$ be the annihilating ideal
of ${\mathcal{F}}$ and ${\Lambda }({\mathcal{F}})$ be the space of the
smooth leafwise forms of ${\mathcal{F}}$. ${\Lambda }({\mathcal{F}})$ is by
definition the quotient space (with the quotient topology) 
\begin{equation}
{\Lambda }^{\ast }(M)\overset{r}{\rightarrow }\frac{{\Lambda }^{\ast }(M)}{I(%
{\mathcal{F}})}={\Lambda }^{\ast }({\mathcal{F}}).  \tag{2.2}
\end{equation}%
Thus the projection $r$ is continuous and \textit{open map\/}. Notice that $%
B^{2}({\mathcal{F}})$ is closed iff 
\begin{equation*}
r^{-1}(B^{2}({\mathcal{F}}))=B^{2}(M)+I^{2}({\mathcal{F}})
\end{equation*}%
is closed. The if direction follows from continuity of $r$. To show the
converse notice that $B^{2}(M)+I^{2}({\mathcal{F}})$ is closed iff 
\begin{equation*}
{\mathcal{A}}={\Lambda }^{2}(M)-(B^{2}(M)+I^{2}({\mathcal{F}}))
\end{equation*}%
is open, thus 
\begin{equation}
r({\mathcal{A}})={\Lambda }^{2}({\mathcal{F}})-B^{2}({\mathcal{F}}) 
\tag{2.3}
\end{equation}%
is open, since $r$ is linear and open. Thus $B^{2}({\mathcal{F}})$ is
closed. Let $D^{2}(M,{\mathcal{F}})=B^{2}(M)\cap I^{2}({\mathcal{F}})$ and
we will show that

\begin{equation}
\Omega =B^{2}(M)+I^{2}({\mathcal{F}})  \tag{2.4}
\end{equation}%
is closed in $\Lambda ^{2}(M)$.For this we show first that $%
I^{2}(\tciFourier )$ is closed in $\Omega .$In fact $\Lambda
^{2}(M)-I^{2}(\tciFourier )$ is open in $\Lambda ^{2}(M)$ since $%
I^{2}(\tciFourier )$

is closed in $\Lambda ^{2}(M)$. Now $\Omega -I^{2}(\tciFourier )=\Omega \cap
(\Lambda ^{2}(M)-I^{2}(\tciFourier ))$ is open in $\Omega $ thus $%
I^{2}(\tciFourier )$ is closed in $\Omega .$Now consider

the restriction $P_{A}^{\ast }:\Omega \rightarrow \Lambda ^{2}(\tciLaplace
^{\ast })$ of $P_{A}$ to $\Omega .$The kernel $\ker P_{A}^{\ast }$ is given
by

\bigskip\ \ \ \ \ \ \ \ \ \ \ \ \ \ \ \ \ \ \ \ \ \ \ \ \ \ \ \ \ \ \ \ \ \
\ \ \ \ \ \ \ \ \ \ \ \ \ \ \ \ \ \ \ \ \ \ \ \ \ \ \ \ \ \ \ \ \ \ \ \ \ \ $%
\ker P_{A}^{\ast }$= $\{\theta =d\alpha +\lambda \Lambda \mu ,P_{A}(\theta
)=P_{A}(d_{f}\alpha )=0\}$\ \ \ \ \ \ \ \ \ \ \ \ \ \ \ \ \ \ \ \ \ \ \ \ \
\ \ \ \ \ \ \ \ \ \ \ \ \ \ \ \ \ \ \ \ \ \ \ \ \ \ \ \ \ \ \ \ \ \ \ \ \ \
\ \ \ \ \ \ \ \ \ \ \ \ \ \ \ \ \ \ \ \ $(2.5$\ $)$

is closed in $\Lambda ^{2}(M).$For if $\theta _{n}\in \ker P_{A}^{\ast }$
converges in the C$^{\infty }$ topology to $\theta \in \Lambda ^{2}(M)$ then 
$\ \theta \in \ker P_{A}^{\ast }.$In fact if $\theta _{n}=d\alpha
_{n}+\lambda _{n}\Lambda \mu $ \ 

where $\alpha _{n}=f_{n}\eta +g_{n}\sigma +h_{n}\mu $ then $P_{A}^{\ast
}(\theta _{n})=$ $P_{A}^{\ast }(d\alpha _{n})=P_{A}^{\ast }(d_{f}\alpha
_{n})=dP_{A}(\alpha _{n})=0$ then $P_{A}(\alpha _{n})=c_{n}\eta $ since $%
d\eta =0$ and $d\sigma =\sigma \Lambda \eta $

in $\ \Lambda (\tciLaplace ^{\ast }).$Now since $\theta _{n}$ converges to $%
\theta $ in the $C^{\infty }$ topology then $P_{A}^{\ast }(\theta _{n})=\
dP_{A}(\alpha _{n})$ $=c_{n}d\eta $ converges in the $C^{\infty }$ topology
to $P_{A}^{\ast }(\theta )=P_{A}(d\alpha )=$

$cd\eta =0$ where $I_{A}(c\eta )=\alpha $ thus $\theta =$ $d\alpha +\lambda
\Lambda \mu $ belongs to $\Omega .$Now I claim that

\bigskip\ \ \ \ \ \ \ \ \ \ \ \ \ \ \ \ \ \ \ \ \ \ \ \ \ \ \ \ \ \ \ \ \ \
\ \ \ \ \ \ \ \ \ \ \ \ \ \ \ \ \ \ \ \ \ \ \ \ \ \ \ \ \ \ \ \ \ \ \ \ \ $%
\Omega =\ker P_{A}^{\ast }+R[d\sigma _{A}]$ \ \ \ \ \ \ \ \ \ \ \ \ \ \ \ \
\ \ \ \ \ \ \ \ \ \ \ \ \ \ \ \ \ \ \ \ \ \ \ \ \ \ \ \ \ \ \ \ \ \ \ \ \ \
\ \ \ \ \ \ \ \ \ \ \ \ \ \ \ \ \ \ \ \ \ \ \ \ \ \ \ \ \ \ \ \ \ \ \ \ \ \
\ \ \ \ \ \ \ \ \ \ \ \ \ \ \ \ \ \ \ \ \ \ \ \ \ \ \ \ \ \ \ \ \ \ \ \ \ \ $%
(2.6)$

For if $\theta \in \Omega $ then $\theta =d\alpha +\lambda \Lambda \mu $
where $\alpha =f\eta +(g+c)\sigma +h\mu $ \ and we may assume that $%
\dint\limits_{M}gd\nu =0$ where $d\nu $ is the measure given by

the canonial volume form $\sigma \Lambda \eta \Lambda \mu .$So $P_{A}^{\ast
}(\theta )=P_{A}^{\ast }(d\alpha )=P_{A}^{\ast }(d_{f}\alpha
)=\{\dint\limits_{M}(Sf-Yg+(g+c))d\nu \}d\sigma =cd\sigma $ so

\ \ \ \ \ \ \ \ \ \ \ \ \ \ \ \ \ \ \ \ \ \ \ \ \ \ \ \ \ \ \ \ \ \ \ \ \ \
\ \ \ \ \ \ \ \ \ \ \ \ \ \ \ \ \ \ \ \ \ $\theta =\theta _{0}+\lambda
\Lambda \mu +cd\sigma _{A}$

where $\theta _{o}=d(f\eta +g\sigma +h\mu )+\lambda \Lambda \mu $ belongs to 
$\ker P_{A}^{\ast }$ \ proving that $\Omega =\ker P_{A}^{\ast }+R[d\alpha
_{A}]$ .Thus $\Omega $ is closed since $\ker P_{A}^{\ast }$ is closed and $%
R[d\sigma _{A}]$

is finite dimensional.

\ \ \ \ \ \ \ \ \ \ \ \ \ \ \ \ \ \ \ \ \ \ \ \ \ \ \ \ \ \ \ \ \ \ \ \ \ \
\ \ \ \ \ \ \ \ \ \ \ \ \ \ \ \ \ \ \ \ \ \ \ \ \ \ \ \ \ \ \ \ \ \ \ \ \ \ 

\bigskip\ \ \ \ \ \ \ \ \ \ \ \ \ \ \ \ \ \ \ \ \ \ \ \ \ \ \ \ \ \ \ \ \ \
\ \ \ \ \ \ \ \ \ \ \ \ \ \ \ \ \ \ \ \ \ \ \ \ \ \ \ \ \ \ \ \ 

\ \ \ \ \ \ \ \ \ \ \ \ \ \ \ \ \ \ \ \ \ \ \ \ \ \ \ \ \ \ \ \ \ \ \ \ \ \
\ \ \ \ \ \ \ \ \ \ \ \ \ \ \ \ \ \ \ \ \ \ \ \ \ \ \ 

\bigskip

\ \ \ \ \ \ \ \ \ \ \ \ \ \ \ \ \ \ \ \ \ \ \ \ \ \ \ \ \ \ \ \ \ \ \ \ \ \
\ \ \ \ \ \ \ \ \ \ \ \ \ \ \ \ \ \ \ \ \ \ \ \ \ \ \ 

\bigskip

\section{The computation of $D^j(M,{\mathcal{F}})$, \quad $2\le j\le 3$}

In [S] we considered the cohomology $H^j(M,{\mathcal{F}})$ of the complex $%
(I({\mathcal{F}}),d)$ and the distinguished closed subspaces $D^j(M,{%
\mathcal{F}})=B^j(M)\cap I({\mathcal{F}})$. From the short exact sequence 
\begin{equation}
0\to I({\mathcal{F}}) \overset{i}{\to} {\Lambda}(M) \overset{r}{\to} {\Lambda%
}({\mathcal{F}})\to0  \tag{3.1}
\end{equation}
we have the long exact cohomology sequence 
\begin{equation*}
0\to H^1(M,R) \overset{r}{\to} H^1({\mathcal{F}}) \overset{\delta}{\to}
D^2(M,{\mathcal{F}})
\end{equation*}
\begin{equation}
\to H^2(M,R) \overset{r}{\to} H^2({\mathcal{F}}) \overset{\delta}{\to} D^3(M,%
{\mathcal{F}})\to 0.  \tag{3.2}
\end{equation}
Now from $H^1({\mathcal{F}})=R[\eta]\oplus H^1(M,R)$ and $H^2({\mathcal{F}}%
)=H^2(M,R)$ and (3.2) we see that 
\begin{equation}
D^2(M,{\mathcal{F}})=R[d\eta] \quad \text{ and }\quad D^3(M,{\mathcal{F}})=0.
\tag{3.3}
\end{equation}

\vskip .2in

\noindent \textbf{Remark.} For any smooth function $f\colon M\rightarrow R$
the 2-form $f\sigma \wedge \eta $ extends to a closed 2-form $w=f\sigma
\wedge \eta +g\sigma \wedge \mu +h\eta \wedge \mu $ i.e. the PDE $\ Uf=Yg-Sh$
has a smooth solution $(g,h)$ for each smooth function $f\colon M\rightarrow
R$. This follows from (3.2) and $D^{3}(M,{\mathcal{F}})=0$.

\pagebreak

\noindent {\Large \textbf{References}}

\begin{itemize}
\item[{[AH]}] J.A. Álvarez and G. Hector, \textit{The dimension of the
leafwise reduced cohomology\/}. Amer. J. Math. \textbf{123} (4) (2001),
604--606.

\item[{[M]}] S. Matsumoto, \textit{Rigidity of locally free Lie group actions
and leafwise cohomology\/}, arXiv: 1002.0393.

\item[{[MM]}] S. Matsumoto and Y. Mitsumatsu, \textit{Leafwise cohomology and
rigidity of certain Lie group actions\/}. Ergod. Th. \& Dynam. Sys. \textbf{%
23} (2003), 61--92.

\item[{[S]}] N.M. dos Santos, \textit{Actions of Lie groups on closed
manifolds\/}. Ergod. Th. \& Dynam. Sys. \textbf{22} (2002), 591--600.
\end{itemize}

\end{document}